\documentclass[journal]{IEEEtran}

\ifCLASSINFOpdf
   \usepackage[pdftex]{graphicx}
   \graphicspath{{./}}
   \DeclareGraphicsExtensions{.pdf,.jpeg,.png,.eps}
\else
\fi

\usepackage[cmex10]{amsmath}
\usepackage{booktabs} 
\usepackage{multirow} 
\usepackage{multicol}
\usepackage{balance}
\usepackage{cite}
\usepackage[switch]{lineno} 
\usepackage{enumitem}
\usepackage{dsfont}
\usepackage{amssymb}
\usepackage{amsthm}

\usepackage[bookmarks=false]{hyperref}
\hypersetup{
  unicode=false,          
  pdftoolbar=true,        
  pdfmenubar=true,        
  pdffitwindow=false,     
  pdfstartview={FitH},    
  pdftitle={My title},    
  pdfauthor={Jun Xing Chin},     
  pdfsubject={Subject},   
  pdfcreator={Creator},   
  pdfproducer={Producer}, 
  pdfkeywords={keyword1} {key2} {key3}, 
  pdfnewwindow=true,      
  colorlinks=true,       
  linkcolor=black,          
  citecolor=black,        
  filecolor=black,      
  urlcolor=black           
}

\usepackage{algorithm}
 \usepackage{algpseudocode}
 
\usepackage[caption=false,font=footnotesize]{subfig}

\begin{document}
%
\title{Machine-Learning Inspired Clustering of Distributed Energy Resources}

\author{
    \IEEEauthorblockN{Jun-Xing~Chin,~\IEEEmembership{Member,~IEEE,}
            and~Gabriela~Hug,~\IEEEmembership{Senior Member,~IEEE}}%

    \thanks{\noindent This work is an outcome of the Future Resilient Systems project at the Singapore-ETH Centre (SEC) supported by the National Research Foundation, Prime Minister's Office, Singapore under its Campus for Research Excellence and Technological Enterprise (CREATE) programme.}

}

\markboth{Preprint}%
{Shell \MakeLowercase{\textit{J. Chin et al.}}: Machine-Learning Inspired Clustering of Distributed Energy Resources}

\maketitle

\begin{abstract}
With the increasing penetration of distributed energy resources (DERs) in distribution grids, their impact on grid operations can no longer be ignored. However, the individual control of these increasingly ubiquitous devices remains a challenge due to their numbers. One solution is to control them in groups via virtual power plants. Previous work has typically focused on optimally dispatching a fixed set of DERs, eschewing methods for optimally determining set membership. Intuitively, these fixed sets may not be optimal from a control perspective. Here, we propose a method to cluster the DERs based on a proxy for their covariances, with the goal of minimising the maximum variance across all DER clusters. This method does not require the enumeration and evaluation of all DER combinations, which are required for brute force techniques. Simulation results show that while there is typically a loss in optimality, it is generally small. More importantly, computational tractability is greatly improved when compared to other methods, which require some form of enumeration and evaluation of the DER combinations.
\end{abstract}


\section{Introduction}
Globally, the penetration of small-scale distributed energy resources (DERs) such as solar photovoltaic panels (PVs), energy storage systems (ESSs), and flexible consumer loads is increasing. Individually, these DERs have negligible impact. However, as they become ubiquitous, a group of these uncontrolled DERs have the potential to negatively impact grid operations, e.g., resulting in over-voltages and increase in operational costs. On the other hand, when properly managed, these DERs have the potential to not only supply energy demand closer to its consumption, but also to enhance system resilience and provide ancillary services for grid operations. Due to their growing numbers, controlling them individually is becoming increasingly challenging. One solution is through aggregating them to form virtual power plants (VPPs), which can then be controlled as aggregated entities to assist in grid operations.

Most works in the literature on VPPs either relate to the real-time dynamic formation of VPPs by agents through an aggregator, e.g. \cite{Yang2018}, or directly through a peer-to-peer process, e.g. \cite{Niesse2014} and \cite{Zhang2021}; or focus on issues pertaining to characterising and leveraging the flexibility offered by these VPPs, e.g. \cite{Mueller2015,Mueller2019}. In \cite{Yang2018}, the authors proposed a scheme where a VPP operator selects and schedules DERs from a pool based on the desired generation profile using stochastic optimisation, while in \cite{Zhang2021} the authors proposed a peer-to-peer framework to form dynamic VPPs based on offers to provide either slow-acting energy based products, or fast-acting power based products to the grid. The framework in \cite{Zhang2021} accounts for an agent's ESS's state-of-charge, and willingness to offer specific services. Note, however, that the term `dynamic' here refers to the commitment of a specific DER to a VPP offer, drawn from an arbitrary pool of DERs.

On the other hand, given a pre-determined set of deterministic DER units, the authors in \cite{Mueller2015} and \cite{Mueller2019} use zonotopes in order to characterise their aggregate flexibility. They then utilise the flexibility values obtained from their method to bid into the power markets, and further propose an algorithm to disaggregate the committed schedules into operational schedules for the individual units \cite{Mueller2019}. Building on that, some of the authors have patented a method to match DERs with zonotopic-flexibility with stochastic DERs \cite{Mueller2020}. This can reduce overall grid uncertainty by having the flexible DERs compensate for the uncertainty of the stochastic resources, at the expense of using said flexibility for assisting grid operations. Also, while intuitively possible, the clustering of DERs based on zonotopic parameters to achieve specific flexibility objectives remain unexplored. More importantly, zonotope-based characterisation methods are not yet readily extendable to stochastic DERs.

Most works on VPPs are based on arbitrary sets of DERs, without discussing how the set memberships are assigned. It is intuitive to see that arbitrary assignments of DERs to VPPs may not lead to the most efficient and optimal groups. In fact, there is a lack of work discussing the clustering/aggregation of DERs, notably when also considering the stochastic nature that most DERs have. Still, clustering problems are not new to power systems researchers, with power grid partitioning and consumer load clustering being well studied problems. Grid partitioning approaches typically relate to topological, structural or electrical distance, with the aim of achieving connected partitions that are highly cohesive internally, but with weak inter-partition relationships (see \cite{Xu2019} for a summary). In \cite{Xu2019}, the authors partition a conventional distribution grid into virtual ``islands'' in order to enhance grid operations and reliability. The authors in \cite{Cotilla-Sanchez2013} use voltage-sensitivity to partition the grid into zones such that intra-zonal transactions have minimal effect on power flows, currents or voltages outside the respective zones. On the other hand, consumer load clustering methods aim to group consumers based on their similarities, in order to improve load forecasting and devise efficient pricing schemes \cite{Zarabie2019,Jia2021}. These methods can either be centralised, e.g., \cite{Zarabie2019}, or distributed, e.g., \cite{Jia2021}, and typically also include demographic and locational data.

Nevertheless, neither grid partitioning nor consumer load clustering works directly result in aggregate profiles that benefit VPP operations and market participation. For VPPs, the DER profiles should ideally be aggregated to result in minimal variance, or maximum flexibility, mimicking those of conventional power plants. In this work, we propose a method to cluster stochastic DERs such that the maximum variance of the resultant aggregates' probability distributions is minimised. The proposed method is based on the  statistical properties of the individual DERs and their correlations to a common external feature. More importantly, this method avoids the need of having to enumerate and evaluate all possible groupings, which does not scale well with the increasing number of DERs.

The rest of the paper is structured as follows: Section \ref{sec:ProbForm} details the problem considered; Section \ref{sec:TheoreticalBounds} presents a possible solution based on theoretical bounds; Section \ref{sec:CovarianceMethod} presents our proposed method; Section \ref{sec:NumEx} presents alternative methods of modelling the optimisation problem, and compares their results against the proposed method using numerical experiments; and Section \ref{sec:Concl} concludes the paper.

\section{Problem Formulation}\label{sec:ProbForm}
Consider the simple problem of clustering a large set of stochastic DERs $\mathcal{N}_{\textit{DER}}$, with the goal of minimising the maximum variance of the resultant profiles across all clusters. Without convoluting the problem with heuristics such as connectivity or topology, and assuming that each DER has to be assigned to a cluster, this can be viewed as a set assignment problem given by the following formulation:
\begin{alignat}{2}
\begin{split}\label{eq:ProbDef}
\underset{x_{i,j}}{\mbox{minimise}} \quad &\underset{j\in\mathcal{N}_{\textit{agg}}}{\max} ~\Gamma\left(\sum_{i\in\mathcal{N}_{\textit{DER}}}  x_{i,j}\right)\\
\mbox{s.t.} \quad & \sum_{j \in \mathcal{N}_{\textit{agg}}} x_{i,j}=1, \quad \quad \forall i\in\mathcal{N}_{\textit{DER}},
\end{split}
\end{alignat}
where $x_{i,j}$ are binary variables assigning DER $i$ to cluster $j$, $\mathcal{N}_{\textit{agg}}$ is the set of desired/possible clusters, and $\Gamma(\cdot)$ gives the variance of the aggregate profile as a function of its constituent DERs. The problem defined by \eqref{eq:ProbDef} can easily be solved by enumerating and evaluating the variances of all possible DER aggregations. However, solving \eqref{eq:ProbDef} using this method quickly becomes computationally intractable as the size of $\mathcal{N}_{\textit{DER}}$ grows. Moreover, the number of DERs are constantly changing, with changes in DER conditions affecting their stochastic properties as well, necessitating the re-clustering of the DERs on a regular basis. Hence, a computationally heavy one-shot solution fails to meet the needs of today's distribution grids, requiring more efficient methods to solving \eqref{eq:ProbDef}. 

\section{Variance Bounds on the DER Aggregates}\label{sec:TheoreticalBounds}
The problem defined by \eqref{eq:ProbDef} can be framed as solving for the variance of the sum of random variables (r.v.s) for each DER aggregate. Assuming that the inter-quantile ranges and medians of the DER aggregates are available, their variances can be estimated using well-known quantile regression methods in order to solve \eqref{eq:ProbDef}. However, this requires the enumeration and evaluation of all potential DER combinations as well. 
 
A potential alternative solution is to first derive bounds on a DER aggregate's quantiles that are based on the quantiles of its member DERs. Then, by using these quantile bounds, further derive bounds on the aggregate's variance. Through minimising them, we ideally also minimise the maximum variance across the DER aggregates. While such a method can indeed be derived, we have found that it is still computationally inefficient. Nevertheless, we think that this carries value for future work, and summarise such an approach in the following. 

In \cite{Watson1986}, Watson and Gordon showed that for specific distribution types, the quantile of a sum of r.v.s can be bounded by a sum of the quantiles of the individual r.v.s. Nonetheless, for the general problem described by \eqref{eq:ProbDef}, the probability distributions of the DERs are not as well-defined, despite there being distribution types that best model specific DER types. Thus, the findings in \cite{Watson1986} are not easily applicable. Fortunately, a potentially viable solution exists in Copula theory, namely through the use of the Fr\'{e}chet-Hoeffding bounds. Let $\{Y_{i}:i\in\mathcal{N}_{\textit{DER}}\}$ be the continuous r.v.s of the DERs, and $Z = \sum_{i\in\mathcal{N}^{k}_\textit{DER}} Y_{i}$ be a new r.v. defined as the sum of the DERs in DER aggregate (subset) $\mathcal{N}^{k}_\textit{DER}$ of size $n_k$. Additionally, let $Q^{a}_{A}$ denote the $a$ quantile of r.v. $A$'s probability distribution and let $F_{A}(\cdot)$ be the cumulative distribution function of $A$. Then, the Fr\'{e}chet-Hoeffding bounds (see \cite{Nelsen2006}, \cite{Schmidt2007} for definition) is given as
\begin{equation}\label{eq:FrechetHoeffding}
     \max\left\{ \sum_{i\in\mathcal{N}^{k}_\textit{DER}} u_{i} + 1 - n_k, 0 \right\} \leq C(\mathbf{u}) \leq \min \left\{u_{1},\cdots,u_{n_k}\right\} .
\end{equation}
$C(\mathbf{u})$ is a copula that represents the dependence structure between the r.v.s $\{Y_{i}: i\in\mathcal{N}_\textit{DER}^{k}\}$, $Q^{u_i}_{Y_i}$ is the $u_i$ quantile of $Y_i$, \emph{i.e.}, $F_{Y_i}(Q^{u_i}_{Y_i}) = u_{i}$, and $\mathbf{u}:=\{u_i: i\in\mathcal{N}_\textit{DER}^{k}\}$. Focusing on the left-hand-side (l.h.s.) of \eqref{eq:FrechetHoeffding}, let $w=\sum_{i\in\mathcal{N}^{k}_\textit{DER}} u_{i} + 1 - n_k$, and fixing $w\ge0$, we have 
\begin{equation} \label{eq:FHlower}
     w \leq C(\mathbf{u}) := \text{Pr}\left(Y_1\leq Q^{u_1}_{Y_1},\cdots, Y_{n_k} \leq Q^{u_{n_k}}_{Y_{n_k}}\right)~,
\end{equation} 
where Pr$(A\leq a)$ is the probability that $A\leq a$. As the probability that the r.v.s $Y_{i}$ are jointly less than their quantiles $Q^{u_i}_{Y_i}$ is at least $w$, the probability that their sum being less than the sum of their quantiles is also bounded from below by $w$,
\begin{alignat}{1}
     w \leq ~\text{Pr}\left(Z \leq \sum_{i\in\mathcal{N}^{k}_\textit{DER}} Q^{u_i}_{Y_i}\right) \Longleftrightarrow Q_{Z}^{w}\leq  \sum_{i\in\mathcal{N}^{k}_\textit{DER}} Q^{u_i}_{Y_i}. \label{eq:sumofquantiles}
\end{alignat}
Equation \eqref{eq:sumofquantiles} is an upper bound on a quantile of $Z$, expressed as the sum of the quantiles of its marginal distributions. 

Using Chebyshev's inequalities, one can derive bounds on a distribution's quantiles, given its mean $\mu$ and standard deviation $\sigma$. For any distribution with a finite mean and non-zero variance, it can be shown\footnote{see \cite{Bagui2004} for details and derivations} that
\begin{alignat*}{1}
\mu_{Z} - \sigma_{Z} \sqrt{\frac{1-w}{w}} &\leq Q^{w}_{Z} \leq \mu_{Z} + \sigma_{Z} \sqrt{\frac{w}{1-w}}~,
\end{alignat*}
which can be rearranged to form lower bounds on its standard deviation (and therefore its variance):
\begin{alignat}{1}\label{eq:stdlower1}
\left(\mu_{Z} - Q^{w}_{Z}\right) \sqrt{\frac{w}{1-w}} &\leq  \sigma_{Z}
\end{alignat}
or
\begin{alignat}{1}\label{eq:stdlower2}
\left(Q^{w}_{Z} - \mu_{Z}\right) \sqrt{\frac{1-w}{w}} &\leq  \sigma_{Z}~.
\end{alignat}
The mean of the aggregate r.v. $Z$, $\mu_{Z}$, can simply be obtained by summing the means of its summands, $\mu_{Z}=\sum_{i\in\mathcal{N}^{k}_\textit{DER}} \mu_{Y_i}$. Depending on whether the non-parametric skew of $Z$ is positive or negative, \emph{i.e.}, whether its mean is greater than or less than its median, either \eqref{eq:stdlower1} or \eqref{eq:stdlower2} will be a tighter bound. As \eqref{eq:sumofquantiles} is an upper bound on $Q_Z$, only \eqref{eq:stdlower1} holds, \emph{i.e.},
\begin{alignat}{1}\label{eq:stdlowerQuantileSum}
\left(\mu_{Z} - \sum_{i\in\mathcal{N}^{k}_\textit{DER}} Q^{u_i}_{Y_i} \right) \sqrt{\frac{w}{1-w}} &\leq  \sigma_{Z}~.
\end{alignat}
While the bound in \eqref{eq:sumofquantiles} can be made tighter by solving for 
\begin{alignat}{2}
\begin{split}\label{eq:minquantilebound}
\underset{w,u_{i}}{\mbox{arg min}} \quad &\sum_{i\in\mathcal{N}^{k}_\textit{DER}} Q^{u_i}_{Y_i}\\
\mbox{s.t.} \quad & w=\sum_{i\in\mathcal{N}^{k}_\textit{DER}} u_{i} + 1 - n_k\\
&w\ge0~,
\end{split}
\end{alignat}
given a particular DER aggregate $\mathcal{N}_\textit{DER}^{k}$, the bound defined by \eqref{eq:stdlowerQuantileSum} may be meaningless if $Z$ has a negative non-parametric skew. Note that the quantile functions in \eqref{eq:minquantilebound} are not always convex \cite{Kibzun2004}. Moreover, in order to derive a solution for \eqref{eq:ProbDef}, the maximum l.h.s. of \eqref{eq:stdlowerQuantileSum} across all clusters needs to be minimised considering all possible DER subset combinations $\mathcal{N}^{k}_\textit{DER} \in \mathcal{P}(\mathcal{N}_\textit{DER})$, where $\mathcal{P}(\mathcal{N}_\textit{DER})$ is the power set of $\mathcal{N}_\textit{DER}$.

Hence, while it is indeed possible to derive a lower bound on the standard deviation (and hence, variance) of a DER aggregate's distribution using the quantiles of its summand r.v.s., there is no guarantee that its actual variance is close to this lower bound. Moreover, in order to obtain the tightest upper bound in \eqref{eq:sumofquantiles}, there is a need to enumerate and solve \eqref{eq:minquantilebound} for all possible DER combinations. Therefore, a bounds-based solution as presented here fails as an efficient approach to solving \eqref{eq:ProbDef}. Thus, motivating a more efficient data-driven method, which shall be presented in the next section.

\section{Variance of a Sum of Random Variables}\label{sec:CovarianceMethod}
From fundamental probability theory, the variance of $Z$, denoted here as $\text{Var}(Z)$ or $\sigma_Z^2$, is given by 
\begin{equation} \label{eq:vardef}
    \text{Var}(Z) = \sum_{i\in\mathcal{N}_\textit{DER}^{k}} \text{Var}(Y_{i}) + 2 \sum_{\substack{\{i,j\} \in \mathcal{N}_\textit{DER}^{k} \\ i<j}} \text{Cov}\left(Y_{i},Y_{j}\right) ~,
\end{equation}
where $\text{Cov}\left(Y_{i},Y_{j}\right)$ is the covariance between r.v.s $Y_i$ and $Y_j$. Directly minimising the maximum of \eqref{eq:vardef} across all clusters would solve \eqref{eq:ProbDef}. However, this would require the use of $C_{2}^{n} := \frac{n !}{2!(n-2)!}$ binary variables and $3 C_{2}^{n}$ linear constraints for each possible DER cluster in the optimisation problem when using conventional mixed-integer linear modelling methods. $C_{2}^{n}$ is the number of 2-combinations for the set of DERs $\mathcal{N}_\textit{DER}$, $n:=|\mathcal{N}_\textit{DER}|$, and $!$ is the factorial operator. For comparison, a brute force approach would require $\sum_{i=1}^{n} C_{i}^{n}$ number of binary variables per desired DER cluster. While the approach using DER covariance is computationally more tractable than a brute force approach, scalability remains a challenge when considering large numbers of DERs.

Using \eqref{eq:vardef} as a starting point, we propose an approach inspired by load-forecasting machine learning techniques. Instead of the covariance, we propose the use of a proxy based on the correlation of the DERs to some external feature. These features should be correlated with all the DERs and could be weather related, e.g., solar irradiance, windspeed, ambient temperature, and humidity. First, we analyse the correlation of these features to the DER profiles, selecting the one that has the highest correlation to all the DER profiles. Ideally, there is a good mix of highly positive and negative correlations amongst the set of DERs to this chosen feature, such that they may cancel out when clustered. Then, we replace the covariance term in \eqref{eq:vardef} with the correlation coefficient multiplied by the variance of the DER distribution. Intuitively, one could also use principal component analysis on the DER data, then dropping all but the first principle component. This component would then be used as the selected feature in the proposed method in order to approximate a solution for \eqref{eq:ProbDef}. Nonetheless, this is left as future work.

A solution for the problem defined in \eqref{eq:ProbDef} can now be approximated by solving the optimisation problem given by:
\begin{alignat}{2}
\begin{split}\label{eq:ClusterMethod}
\underset{x_{i,j}}{\mbox{minimise}} \quad &\underset{j\in\mathcal{N}_{\textit{agg}}}{\max} ~  \sum_{i\in\mathcal{N}_\textit{DER}} ~x_{i,j} \sigma_{Y_{i}}^2 + \\
& \quad \quad ~ \underbrace{ \left| \sum_{i\in\mathcal{N}_\textit{DER}}x_{i,j} \text{Corr}(Y_{i})\sigma_{Y_{i}}^2 \right|}_{\text{proxy for covariance}} ~\\
\mbox{s.t.} \quad & \sum_{j \in \mathcal{N}_{\textit{agg}}} x_{i,j}=1, \quad \quad \forall i\in\mathcal{N}_{\textit{DER}},
\end{split}
\end{alignat}
where $\text{Corr}(Y_{i})$ is the correlation of DER $Y_i$ to the selected external feature, and $\sigma_{Y_i}^2$ is the variance of DER $i$. In \eqref{eq:vardef}, the lowest variance for $Z$ is achieved when the sum of the covariance terms are most negative. However, when using the proposed covariance proxy, there is a need to take its absolute value, as a highly negative term does not necessarily equate to the lowest covariance between the DERs. More specifically, if all the DERs within a cluster are highly negatively correlated, it actually equates to their covariances being highly positive instead, which is to be avoided. Hence, the goal here is to minimise the absolute covariance proxy term when summed across the summand DERs, either through a combination of DERs with matching positive and negative correlations to the external feature that cancel out, or a combination of DERs that result in the minimal covariance proxy terms across all clusters. Note that solving \eqref{eq:ClusterMethod} does not guarantee an optimal solution for \eqref{eq:ProbDef}, but should outperform most random cluster assignments. While bounds do exist on the correlation between two variables given their correlation to a third common variable (feature), these bounds do not readily extend to multiple variables, as required in this problem.

As the problem in \eqref{eq:ClusterMethod} cannot be directly solved by an optimisation solver, it is reformulated as the following:
\begin{alignat}{2}\label{eq:ClusterOpt}
\underset{x_{i,j},y,z}{\mbox{minimise}} \quad & ay + bz\\
\mbox{s.t.} \quad & \sum_{j \in \mathcal{N}_{\textit{agg}}} x_{i,j}=1, \quad \quad &&\forall i\in\mathcal{N}_{\textit{DER}}\notag\\
& \sum_{i\in\mathcal{N}_\textit{DER}} x_{i,j} \sigma_{Y_{i}}^2 \leq y,  \quad \quad &&\forall j\in\mathcal{N}_{\textit{agg}} \notag\\
& \sum_{i\in\mathcal{N}_\textit{DER}} x_{i,j} \text{Corr}(Y_{i})\sigma_{Y_{i}}^2 \leq z, \quad \quad &&\forall j\in\mathcal{N}_{\textit{agg}}\notag\\
& \sum_{i\in\mathcal{N}_\textit{DER}} x_{i,j} \text{Corr}(Y_{i})\sigma_{Y_{i}}^2 \geq -z, \quad \quad &&\forall j\in\mathcal{N}_{\textit{agg}}\notag,
\end{alignat}
where $a$ and $b$ are weights that can be used to prioritise either the variance or covariance proxy terms, and $y$ and $z$ are dummy variables used to reformulate the minimax optimisation problem. As stated previously, a highly negative covariance proxy term does not necessarily result in a highly negative covariance, but would result in a lower objective function value. As such, it is necessary to minimise these terms separately, requiring separate dummy variables in the optimisation problem.

Based on the Cauchy-Schwarz' inequality, the covariance between two r.v.s is bounded by
\begin{equation*}
   \left| \text{Cov}\left(Y_{i},Y_{j}\right) \right|  \leq \sqrt{\sigma^2_{Y_{i}} \sigma^2_{Y_{j}}}~.
\end{equation*}
Hence, the covariance would be of similar magnitude as the variance of the individual r.v.s. Given that the correlation of the DERs to the external feature is within the range of $[-1,1]$, the coefficients $a$ and $b$ can generally, both be set to one. 

Note that \eqref{eq:ClusterOpt} only enforces a maximum number of clusters, but not a minimum; and that while \eqref{eq:ClusterOpt} still contains binary variables, their quantity is limited to the number of DERs (size of set $\mathcal{N}_\textit{DER}$), making it much more computationally tractable when compared to a brute force approach or an optimisation problem using \eqref{eq:vardef}.

\section{Case Studies and Numerical Experiments}\label{sec:NumEx}
To the best of our knowledge, no traditional methods for solving this problem exists. Hence, to validate our proposed clustering approach, we compare it to brute force and variance formula based methods using numerical simulations. For ease of comparison, the results of the three optimisation based clustering methods are benchmarked against results from Monte Carlo  simulations (MC sims) that randomly assign the DER clusters. The modelling of the brute force and the variance formula based optimisation problems are summarised in the next two sub-sections, followed by numerical experiments and their discussion.

\subsection{Brute Force Method of Clustering the DERs}
A brute force method for clustering the DERs involves enumerating and evaluating the variances for all possible combinations of the DERs. These combinations and the associated cluster variances can then be formulated as a combinatorial optimisation problem. Define $\sigma_{Z_k}^{2}$ as the variance of the sum of the r.v.s in the $k$-th DER subset, which is estimated from available data. Further, let $c_{k,j}$ be binary variables that assign the $k$-th DER subset to cluster $j$, $\mathcal{C}_{k} \in \mathcal{P}(\mathcal{N}_\textit{DER})$ the $k$-th DER subset, and $N_k := |\mathcal{P}(\mathcal{N}_\textit{DER}|)$ the number of possible DER subsets. Then, a straightforward brute force based optimisation problem (\textit{brute force model}) can be formulated as:
\begin{alignat}{2}\label{eq:BruteOpt}
\underset{x_{i,j},c_{k,j},z}{\mbox{minimise}} \quad & z\\
\mbox{s.t.} \quad & \sum_{j \in \mathcal{N}_{\textit{agg}}} x_{i,j}=1, ~ &&\forall i\in\mathcal{N}_{\textit{DER}}\\
& \sum_{k=1}^{N_k} c_{k,j} \sigma_{Z_k}^{2}  \leq z, ~ &&\forall j\in\mathcal{N}_{\textit{agg}}\\
& \sum_{i \in \mathcal{C}_{k}} x_{i,j} \ge c_{k,j}|\mathcal{C}_k|,~ &&\forall j\in\mathcal{N}_{\textit{agg}}, \label{eq:combomember}\\ 
& && ~\,k=1,\cdots,N_k \notag\\
& \sum_{k=1}^{N_k} c_{k,j} |\mathcal{C}_k| \ge \sum_{i\in\mathcal{N}_\textit{DER}} x_{i,j}, ~ &&\forall j\in\mathcal{N}_{\textit{agg}}\label{eq:correctcombo}\\
& \sum_{k=1}^{N_k} c_{k,j}  \leq 1, ~ &&\forall j\in\mathcal{N}_{\textit{agg}}\label{eq:maxcombo}.
\end{alignat}
Here, constraints \eqref{eq:combomember} to \eqref{eq:maxcombo} ensure that the correct variance values are assigned to a cluster given the DERs assigned to it.

\subsection{Clustering of DERs based on Multi-variable Variance}
For this approach, which we denote as the \textit{covariance model}, \eqref{eq:vardef} is used in the objective function to calculate the variance for each DER cluster. Let $\mathcal{N}_\textit{DER}^{[2]}$ be the set of 2-subsets of $\mathcal{N}_\textit{DER}$, with size $N_c:=\left|\mathcal{N}_\textit{DER}^{[2]}\right|$, and define $c_{k,j}$ as binary variables that indicate whether the DER pair in the 2-subset $\mathcal{C}_k \in \mathcal{N}_\textit{DER}^{[2]}$ is in cluster $j$. Additionally, let $\text{Cov}\left(\mathcal{C}_{k}\right)$ be the covariance between the DER pair in subset $\mathcal{C}_k$, which is estimated beforehand and entered into the optimisation problem as constants. Then an optimisation problem that solves \eqref{eq:ProbDef}, which is based directly on \eqref{eq:vardef} is given by:

\begin{alignat}{2}\label{eq:CovOpt}
\underset{x_{i,j},c_{k,j},z}{\mbox{minimise}} \quad & z\\
\mbox{s.t.} \quad & \sum_{j \in \mathcal{N}_{\textit{agg}}} x_{i,j}=1, ~ &&\forall i\in\mathcal{N}_{\textit{DER}}\\
& \sum_{i\in\mathcal{N}_\textit{DER}} x_{i,j}\sigma_{Y_i}^2 + 2 \sum_{k=1}^{N_c} c_{k,j} &&\text{Cov}\left(\mathcal{C}_{k}\right) \leq z, \notag\\
& ~ &&\forall j\in\mathcal{N}_{\textit{agg}} \label{eq:dummyVarCov}\\
& \sum_{i \in \mathcal{C}_{k}} x_{i,j} \ge 2 ~c_{k,j}, ~ &&\forall j\in\mathcal{N}_{\textit{agg}},\notag\\ 
& && ~\,k=1,\cdots,N_c \label{eq:pairmemberUp}\\
&  c_{k,j}   \ge \sum_{i\in\mathcal{C}_{k}} x_{i,j} - 1, ~ &&\forall j\in\mathcal{N}_{\textit{agg}},\notag\\ 
& && ~\,k=1,\cdots,N_c \label{eq:pairmemberDown}.
\end{alignat}
The constraints \eqref{eq:pairmemberUp} and \eqref{eq:pairmemberDown} ensure that the covariance of the DER pair in $\mathcal{C}_{k}$ are included correctly in \eqref{eq:dummyVarCov} if both DERs in the pair are in cluster $j$.

\subsection{Numerical Experiments}
For the numerical experiments, we use 15-minute building-level power consumption and PV generation data from Elektrizit\"{a}tswerke des Kantons Z\"{u}rich (EKZ), a Swiss distribution system operator. There are a total of $14$ PV generators and $36$ building-level consumption profiles in this dataset. For the external feature data, local weather data (global radiation, windspeed, temperature, and relative humidity) obtained from the Swiss Federal Office of Meteorology and Climatology (MeteoSwiss) is used. It is assumed in this study that the buildings in this dataset are participants in a demand response (DR) scheme with flexible loads, and the goal is to group these DR participants and PV generators to obtain DER clusters with minimal variance. The time-of-day and the seasons greatly impact PV production and consumer load behaviour, and thus, the cluster assignments. Here, we consider the clustering problem for the spring and summer months (31 March to 27 October), between the hours of 09:00 and 18:00.  

\begin{figure}
    \centering
    \includegraphics[trim=0cm 0.5cm 1.5cm 2cm, clip=true, width=0.90\columnwidth]{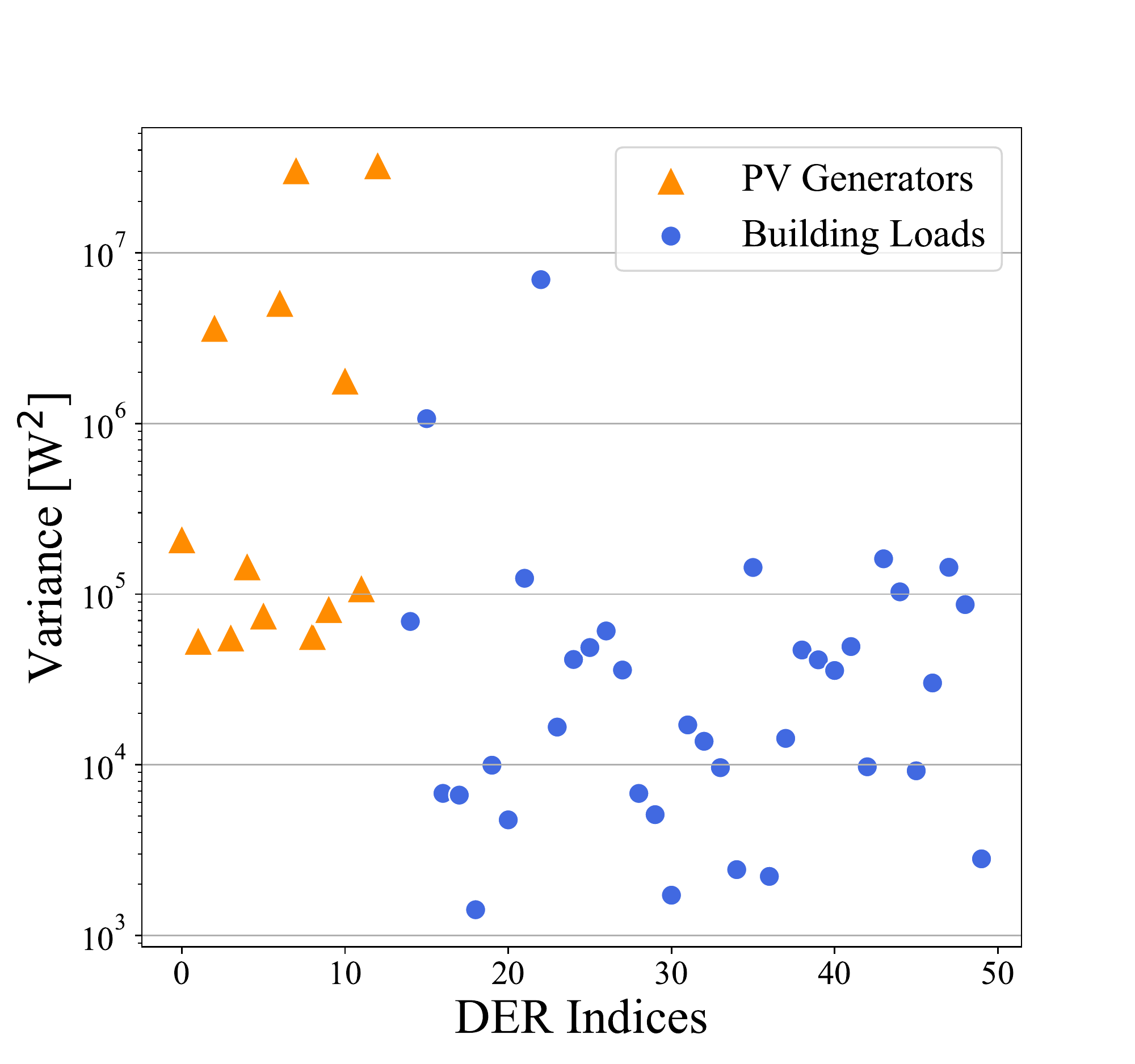}
    \vspace{-0.1cm}
    \caption{Variance of DERs in the EKZ Dataset}
    \label{fig:DERVar}
\end{figure}

\begin{figure}
    \centering
    \includegraphics[trim=0cm 0.5cm 1.5cm 1cm, width=0.90\columnwidth]{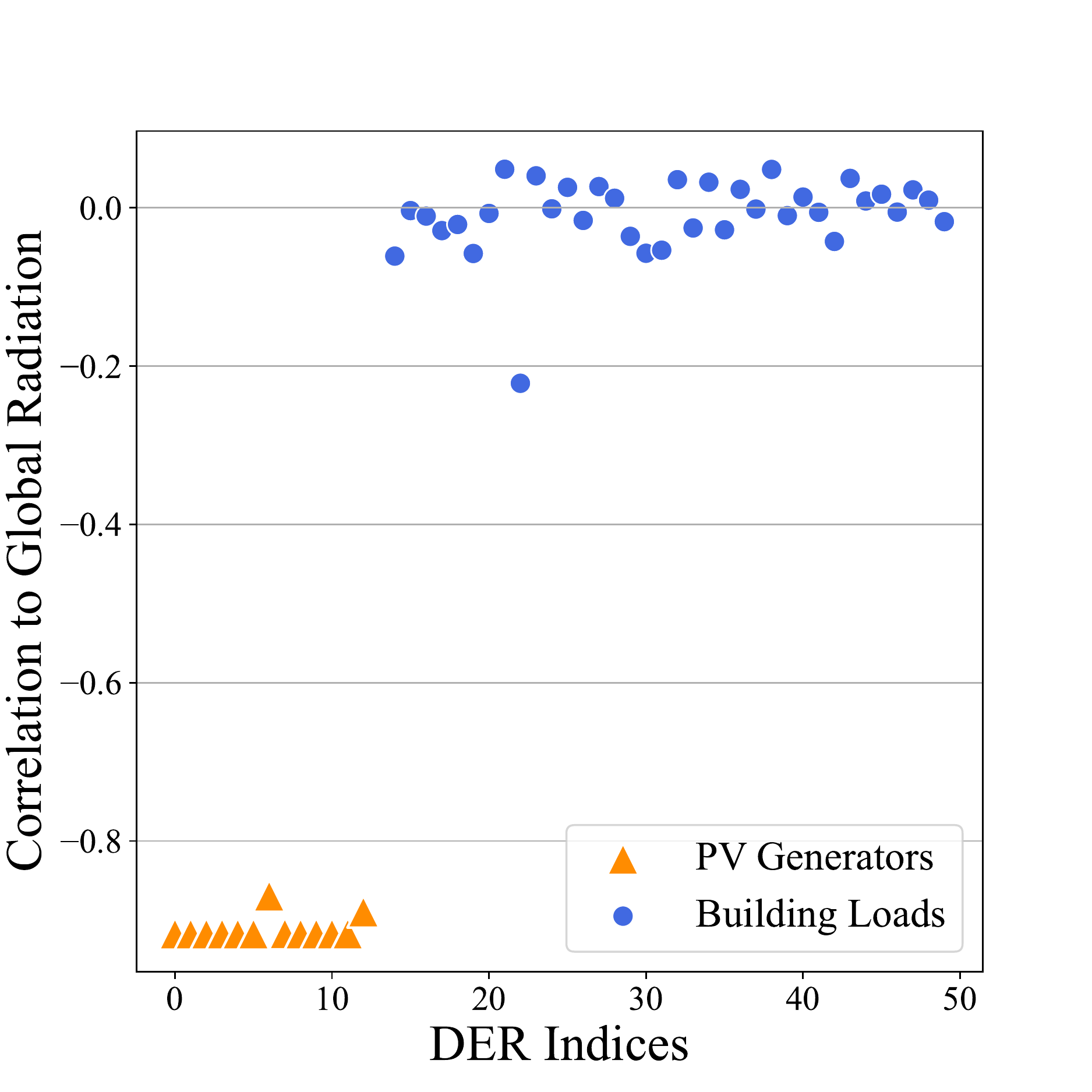} \vspace{-0.1cm}
    \caption{Correlation of DERs to Global Radiation}
    \label{fig:DERCorr}
    \vspace{-0.3cm}
\end{figure}

\begin{figure}
    \centering
    \includegraphics[trim=1cm 0.5cm 2.85cm 2cm, clip=true, width=0.98\columnwidth]{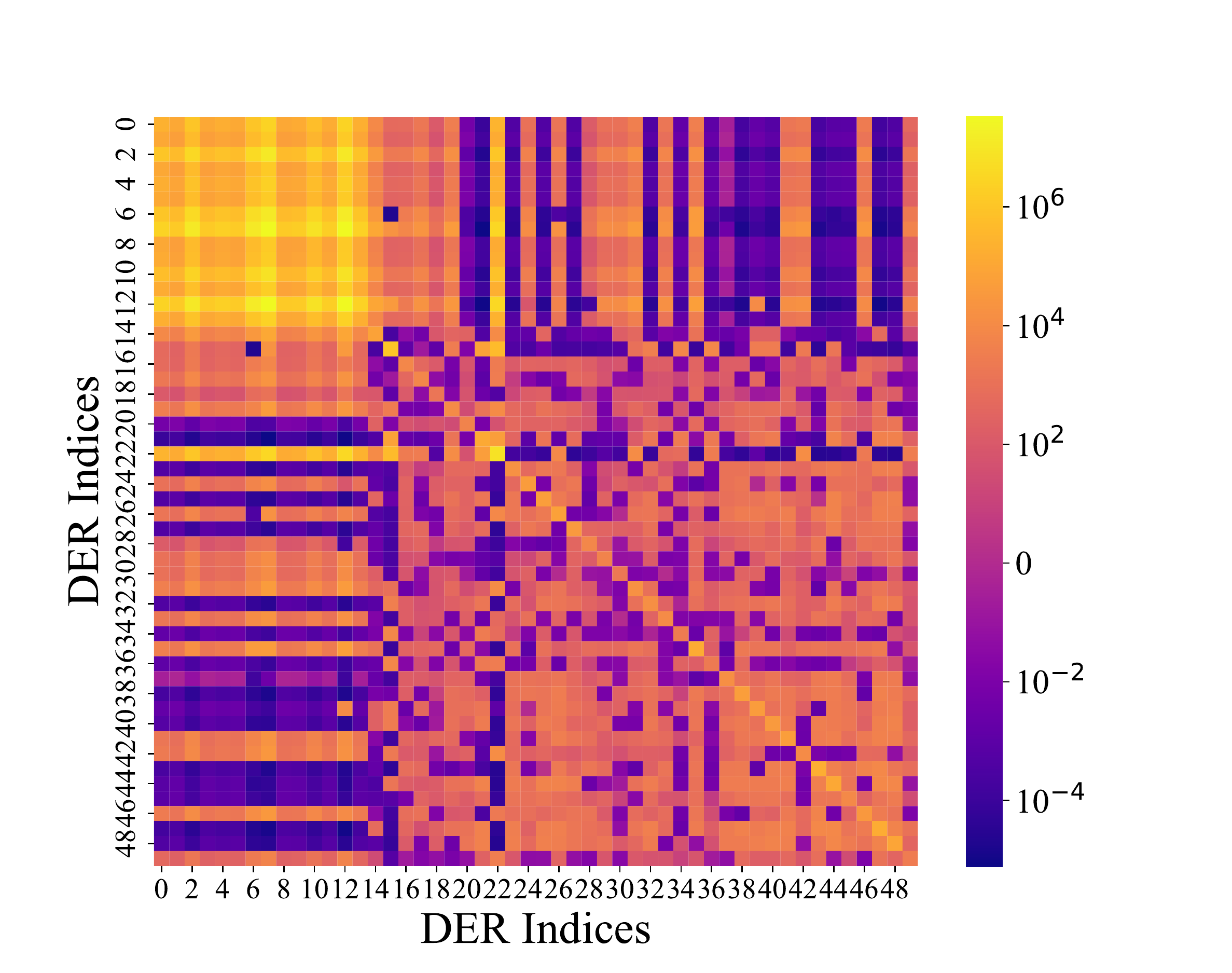} \vspace{-0.1cm}
    \caption{Covariance of DERs in W$^{2}$}
    \label{fig:DERCovar}
    \vspace{-0.3cm}
\end{figure}

Based on a quick correlation analysis, it was found that the global (solar) radiation had the most influence on the PV and load distributions. Hence, the correlations between the DERs and global radiation were used as the input for the proposed optimisation method in this case study. Fig. \ref{fig:DERVar} shows the variance for the individual DERs in the dataset used, while Fig. \ref{fig:DERCorr} illustrates the correlation of the DERs to the global radiation values recorded at a nearby weather station. As can be seen in Fig. \ref{fig:DERCorr}, the PV generators are highly correlated with global radiation; values are negative as generation is taken as negative from the grid perspective. Most of the consumer loads have low correlations with the global radiation, indicating that the buildings are probably well insulated against solar heat gains, and do not utilise energy intensive cooling appliances. Interestingly, some consumer load is negatively correlated with the global radiation, which could be caused by a reduction in electric heating. The covariances among the DERs are summarised in Fig. \ref{fig:DERCovar} as a heat map. As can be seen in Fig. \ref{fig:DERCovar}, the PV generators have high positive covariances, but have in turn, negative covariances with most of the consumer loads.

\begin{table}
\renewcommand{\arraystretch}{1.1}
\centering
\caption{Parameters and Setup for Numerical Experiments}
\label{tab:SimSetup}
\begin{tabular}{|c|c|} \hline 
No. of PV DERs per run          & $8$       \\ \hline
No. of load DERs per run        & $8$       \\ \hline
Max. no. of clusters            & $4$       \\ \hline
No. of random runs, $m_\textit{opt}$   & $250$     \\ \hline
No. of MC sims/run, $n_\textit{mc}$    & $100,000$  \\ \hline
Coefficients ${a,b}$            & ${1,1}$   \\ \hline
\end{tabular}
\end{table}

In order to validate that the results obtained were not specific to a particular set of DERs, the optimisation problems were repeated for $m_\textit{opt}$ times using $16$ DERs randomly drawn from the $50$ that are present in the dataset. Unless otherwise stated, the basic parameters and setup used are as shown in Table \ref{tab:SimSetup}. Using $n_\textit{mc}=100,000$ MC simulations, a distribution of the maximum variance from randomly assigned DER clusters was constructed. Choosing a large number of $n_\textit{mc}$ ensures that a fairly accurate assessment of randomly assigned clusters is achieved, as the uncertainty in the cumulative probability values is proportional to $n_\textit{mc}^{-1/2}$ \cite{Lerche2005}. The performance of the three optimisation methods were then assessed by gauging how they fared against the randomly assigned clusters, \emph{i.e.}, the percentile of MC simulation runs that have a lower maximum variance than the optimisation results.

The computational tractability of optimisation methods were also assessed by measuring the solution time using the Gurobi 9.1.2 solver in Python; running on dedicated Windows Server 2019 Virtual Machines, each with eighteen 3.1 GHz Intel Xeon Gold cores and 256GB of RAM. The enumeration and evaluation of all possible DER combinations for the brute force and covariance methods were not considered in the computational tractability assessment because such operations could potentially be made faster by more efficient implementations. 

\begin{figure}
    \centering
    \includegraphics[trim=0cm 0.5cm 1.5cm 2cm, clip=true, width=0.90\columnwidth]{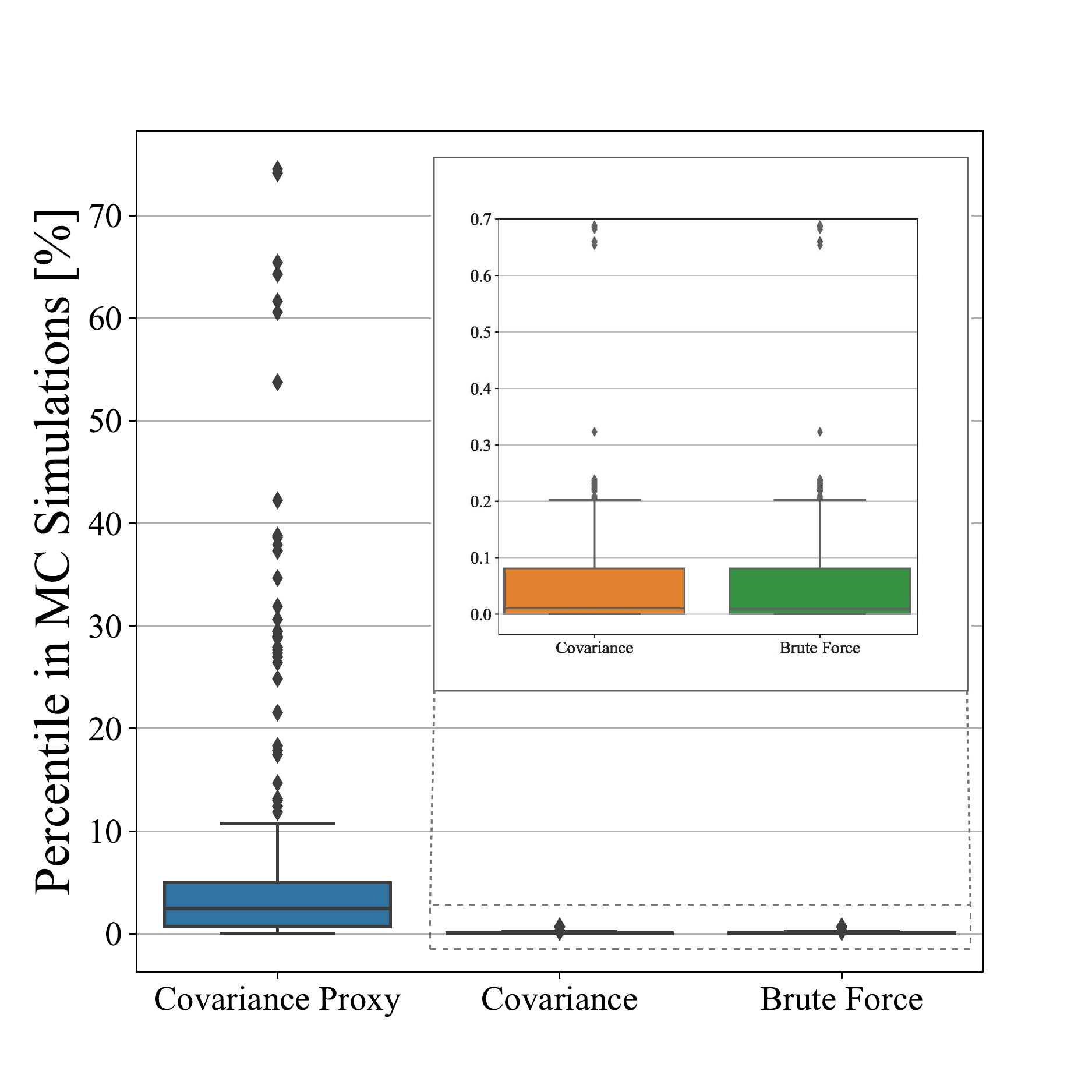} 
    \vspace{-0.1cm}
    \caption{Maximum variance of the optimisation solutions as percentiles of the Monte Carlo simulation results}
    \label{fig:ClusterVar}
\end{figure}

\begin{figure}
    \centering
    \includegraphics[trim=0cm 0.5cm 1.5cm 2cm, clip=true, width=0.90\columnwidth]{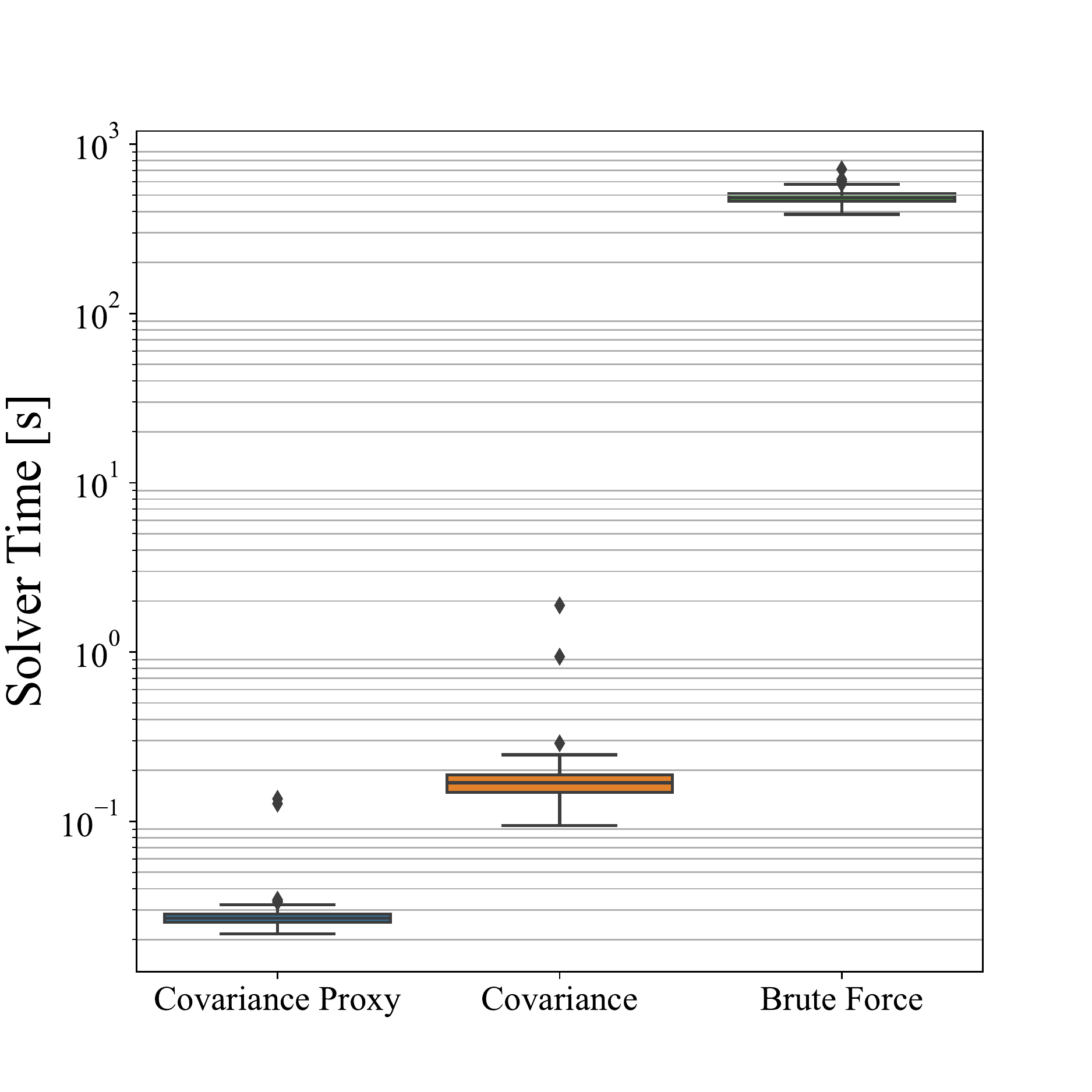} 
    \vspace{-0.1cm}
    \caption{Solution time for the optimisation methods}
    \label{fig:Runtime}
    \vspace{-0.3cm}
\end{figure}

Fig. \ref{fig:ClusterVar} and \ref{fig:Runtime} illustrate the percentiles of the MC simulation on which the solutions from the optimisation methods fall, and the time required to solve, respectively. A lower value for the former means lower cluster variance, which is desired.
As illustrated, it can be seen that the proposed method outperforms at least $50\%$ of randomly assigned clusters in most cases ($97.2\%$ or $243$ out of $250$ runs). On average, it outperformed or was equal to $93.02\%$ of the random MC cluster assignments. While there can be instances where the proposed method fares worse than most of the randomly assigned clusters, these instances are few in our experiments. Also, as expected, the covariance and brute force models are on par with the most optimal clustering results from the MC cluster assignments, as seen in Fig. \ref{fig:ClusterVar}. The fact that they are not always at the zeroth percentile in all cases is because there are multiple optimal MC cluster assignments. Note, however, that the covariance model may exhibit slight sub-optimality when compared to the brute force model in some cases. This occurs when there are issues with data completeness, which could also affect the proposed covariance proxy method as well. When assessing the DER covariances and cluster variances, the availability of data is limited by the least complete DER time series. Then, it follows that when issues with data completeness arise, the cluster variances would be estimated based on the least amount of data, followed be the DER covariances and the individual DER variances. This could result in discrepancies, which objectively does not necessarily indicate that the brute force model outperforms the covariance model, but could affect the results in a numerical study. Also, note that even if a resultant cluster profile has zero variance, the resultant mean is not necessarily zero.

More importantly, as can be seen in Fig. \ref{fig:Runtime}, computational tractability is greatly improved with the proposed method, despite the clustering problem only considering $16$ DERs. The brute force model required, on average, $488$ seconds to solve, while the proposed method only needed $0.028$ seconds. The brute force model also had a maximum solution time of $711$ seconds despite the few DERs. While the covariance model had an average solution time of $0.179$ seconds here, this increases significantly with the increase in the number of DERs. To study this, we included additional DERs and possible clusters in the optimisation problems. The brute force model was not scaled beyond $21$ DERs due to computational tractability issues. With the exception of the brute force model, these experiments were repeated $350$ times each, with different sets of randomly selected DERs. Performance-wise, the proposed method exhibited results comparable to that of the $16$ DER case for the scaled-up clustering problems. Table \ref{tab:scalability} summarises the maximum solution times required.

\begin{table}[!t]
\renewcommand{\arraystretch}{1}
\centering
\caption{Scalability of the three optimisation methods}
\label{tab:scalability}
\begin{tabular}{|l|c|r|r|} \hline
~                           &\multicolumn{3}{c|}{Max Opt. Solution Time [seconds]}\\ \hline
DERs and Cluster Quantity   & Covar. Proxy      & Covar.        & Brute Force   \\ \hline
10 PV, 10 loads; 4 Clust.   & $0.0312$          & $2.683$      & $171,052$     \\ \hline
11 PV, 10 loads; 4 Clust.   & $0.0382$          & $5.597$      & $1,490,401$   \\ \hline
12 PV, 12 loads; 4 Clust.   & $0.0473$          & $9.688$      & -           \\ \hline
14 PV, 14 loads; 4 Clust.   & $0.0467$          & $0.9558$      & -           \\ \hline
14 PV, 14 loads; 8 Clust    & $0.7836$          & $30.56$       & -           \\ \hline
14 PV, 14 loads; 16 Clust   & $0.1367$          & $16.24$       & -           \\ \hline
14 PV, 14 loads; 24 Clust.  & $0.4850$          & $2,023$       & -           \\ \hline
14 PV, 20 loads; 24 Clust.  & $0.3129$          & $4,001$       & -           \\ \hline
14 PV, 26 loads; 24 Clust.  & $0.4915$          & $25,733$      & -           \\ \hline
\end{tabular}
\vspace{-0.5cm}
\end{table}

The brute force model required approximately $47.5$ hours to solve with $20$ DERs and $4$ possible clusters. This increases exponentially to a maximum of $414$ hours ($17.25$ days) when considering $21$ DERs and $4$ possible clusters, which makes scaling this model computationally challenging. While the covariance model remains tractable even up to $40$ DERs and $24$ clusters, requiring approximately $7$ hours to solve, the solution time has increased exponentially when compared to a setting with just $20$ DERs and $4$ clusters. It can reasonably be expected that further scaling the covariance model would be computationally challenging. On the other hand, the proposed covariance proxy model only required marginally more time to solve, remaining highly tractable, and thus, scalable. Moreover, if one were to also consider the times required to enumerate and evaluate the DER combinations for the brute force and covariance models, this difference in computational tractability is further exacerbated. Note that solution time does not increase monotonically with small increases in problem size due to the mixed-integer nature of the problem.

While there is generally a loss of optimality for the proposed covariance proxy model, this is made up for by its significant improvement in computational tractability and scalability. In the near future, distribution grids are expected to host hundreds, if not thousands, of DERs. Hence, among the three models compared here, only the covariance  proxy model has been shown to be viable for clustering these numerous DERs. More importantly, the tractability of the proposed method allows for real-time reclustering of the DERs, enabling dynamic updates that incorporate operational changes.

\section{Conclusion and Outlook}\label{sec:Concl}
In this work, we have shown that by using a proxy for the covariance of the DERs, one can significantly improve the computational tractability of clustering stochastic DERs based on their probability distributions to minimise the maximum variance of the resultant DER aggregate profiles. Despite the loss in optimality and the lack of guarantees that it will always outperform randomly assigned clusters, the proposed method is generally on par with brute force and variance formula based cluster assignment methods. Hence, the proposed covariance proxy based model is a viable and practical approach in clustering the increasingly ubiquitous DERs in distribution grids. Future work will focus on achieving other stochastic properties of the DER aggregates, developing methods to improve or bound the worst case performance of the proposed method, and potentially clustering the DERs based on specific flexibility requirements using their zonotopic parameters.

\section*{Acknowledgements}
The authors would like to thank Elektrizit\"{a}tswerke des Kantons Z\"{u}rich (EKZ) for providing the distribution grid data that forms the basis of the numerical analayses in this work. 

\bibliographystyle{IEEEtran}
\bibliography{./DERCluster}

\end{document}